\documentclass[11pt]{article}
\usepackage{setspace}
\usepackage{textcomp}
\usepackage{amsmath, amsthm, amssymb}
\usepackage{hyperref}
\usepackage{latexsym}
\usepackage{amssymb}
\usepackage{mathrsfs}
\usepackage{rawfonts}
\input{prepictex}
\input{pictex}
\input{postpictex}

\newtheorem{thm}{Theorem}[section]
\newtheorem{lem}[thm]{Lemma}
\newtheorem{prop}[thm]{Proposition}

\def\ad{\mathrm{ad}}
\def\bct{\begin{center}}
\def\ect{\end{center}}
\def\beg{\begin}

\def\<{\langle}
\def\>{\rangle}
\def\lra{\longrightarrow}
\def\mbb{\mathbb}
\def\mbbp{\mathbb P}
\def\mc{\mathcal}
\def\mco{\mathcal O}
\def\mcz{\mathcal Z}
\def\mf{\mathfrak}
\def\mrm{\mathrm}

\def\tn{\textnormal}
\def\wt{\widetilde}
\def\wtz{\wt{Z}}
\newtheorem{rmk}[thm]{Remark}

\def\bcp{{\Bbb C\Bbb P}}

\newcounter{exmp}

\title
{Scalar Curvature and Connected Sums of Self-Dual 4-Manifolds}

\author{Mustafa Kalafat}

\begin{document}
\maketitle

\begin{abstract} In \cite{DF} under a reasonable vanishing hypothesis, Donaldson
and Friedman proved that the connected sum of two self-dual
Riemannian 4-Manifolds is again self-dual. Here we prove that the
same result can be extended over to the positive scalar curvature
case. 

The proof is based on the twistor theory. First a technical
vanishing theorem is established by using an appropriate spectral
sequence, then Green's functions and Serre-Horrocks vector bundle
constructions are used to detect the sign of the scalar curvature.

\end{abstract}


\section{Introduction}



\index{Weyl Curvature Tensor} Let $(M,g)$ be an oriented Riemannian
n-manifold. Then by raising an index, the Riemann curvature tensor
at any point can be viewed as an operator $\mathcal R : \Lambda^2M
\to \Lambda^2M$ hence an element of $S^2\Lambda^2M$. It satisfies
the algebraic Bianchi identity hence lies in the vector space of
{\em algebraic curvature tensors}. This space is an $O(n)$-module
and has an orthogonal decomposition into irreducible subspaces for
$n\geq 4$. Accordingly the Riemann curvature operator decomposes as:
$$\mathcal R = U \oplus Z \oplus W$$
where
$$U={s\over 2n(n-1)}g\bullet g \hspace{.5cm}\textnormal{and} \hspace{.5cm} Z={1\over n-2}\stackrel{\circ}{Ric}\bullet g$$ $s$ is
the scalar curvature, $\stackrel{\circ}{Ric}=Ric-{s\over n}g$ is the
trace-free Ricci tensor, "$\bullet$" is the Kulkarni-Nomizu product,
and $W$ is the {\em Weyl Tensor} which is defined to be what is left
over from the first two pieces.

When we restrict ourselves to dimension $n=4$, the Hodge Star
operator $* : \Lambda^2 \to \Lambda^2$ is an involution and has $\pm
1$-eigenspaces decomposing the space of two forms as $\Lambda^2 =
\Lambda_+^2\oplus\Lambda_-^2$, yielding a decomposition of any
operator acting on this space. In particular $W_\pm : \Lambda_\pm^2
\to \Lambda_\pm^2$ is called the self-dual and anti-self-dual pieces
of the Weyl curvature operator. And we call $g$ to be {\em
self-dual}(resp. {\em anti-self-dual}) {\em metric} if $W_-$(resp.
$W_+$) vanishes. In this case \cite{ahs} construct a complex
$3$-manifold $Z$ called the {\em Twistor Space} of $(M^4,g)$, which
comes with a fibration by holomorphically embedded rational curves
:\vspace{.2cm}

\hspace{3.5cm} \begin{tabular}{cccc}
$\mathbb{CP}_1$ & $\to$ & $Z$ & Complex $3$-manifold \\
&&$\downarrow$& \\
&&$M^4$ & Riemannian $4$-manifold
\end{tabular}
\vspace{.2cm}

This construction drew the attention of geometers, and many examples
of Self-Dual metrics and related Twistor spaces were given
afterwards. One result proved to be a quite effective way to produce
infinitely many examples and became a cornerstone in the field :\\\\
{\bf Theorem \ref{thmb}} (Donaldson-Friedman,1989,\cite{DF}). {\em
If $(M_1,g_1)$ and $(M_2,g_2)$ are compact self-dual Riemannian
4-manifolds with $H^2(Z_i,\mco(TZ_i))=0$

Then~~~~$M_1\#M_2$ also admits a self-dual metric.}\\

The idea of the proof is to work upstairs in the complex category
rather than downstairs. One glues the blown up twistor spaces from
their exceptional divisors to obtain a singular complex space
$Z_0=\wt Z_1 \cup_Q \wt Z_2$. Then using the Kodaira-Spencer
deformation theory extended by R.Friedman to singular spaces, one
obtains a smooth complex manifold, which turns out to be the twistor
space of the connected sum.

When working in differential geometry, one often deals with the
moduli space of certain kind of metrics. The situation is also the
same for the self-dual theory. Many people obtained results on the
space of positive scalar curvature self-dual(PSC-SD) metrics on
various kinds of manifolds. Since the positivity of the scalar
curvature imposes some topological restrictions on the moduli space,
people often find it convenient to work under this assumption.

However one realizes that there is no connected sum theorem for
self-dual positive scalar curvature metrics. Donaldson-Friedman
Theorem(\ref{thmb}) does not make any statement about the scalar
curvature of the metrics produced. Therefore we attacked the problem
of determining the sign of the scalar curvature for the metrics
produced over the connected sum, beginning by proving the
following, using the techniques similar to that of \cite{lom}:\\\\
{\bf Theorem \ref{vanishing}} (Vanishing Theorem). {\em Let $\omega
: {\mathcal Z} \to {\mathcal U}$ be a $1$-parameter standard
deformation of $Z_0$, where $Z_0$ is as in Theorem (\ref{thmb}), and
${\mathcal U}\subset \mathbb C$ is a neighborhood of the origin. Let
$L\to {\mathcal Z}$ be the holomorphic line bundle defined by
$${\mathcal O}(L^*)
={\mathscr I}_{\wt{Z}_1}(K_{\mathcal Z}^{1/2}).$$ If $(M_i,[g_i])$
has positive scalar curvature, then  by possibly replacing
${\mathcal U}$ with a smaller
 neighborhood  of $0\in \mathbb C$ and simultaneously replacing
 ${\mathcal Z} $ with its inverse image, we can arrange for
 our complex $4$-fold $\mathcal Z$ to satisfy
 $$H^1 ({\mathcal Z} , {\mathcal O}(L^*))=H^2 ({\mathcal Z} , {\mathcal
   O}(L^*))=0.$$}

The proof makes use of the Leray Spectral Sequence, homological
algebra and Kodaira-Spencer deformation theory, involving many
steps. Using this technical theorem next we prove that the
Donaldson-Friedman Theorem can be generalized to the
positive scalar curvature(PSC) case :\\\\
{\bf Theorem \ref{pos}.} {\em Let $(M_1,g_1)$ and $(M_2,g_2)$be
compact self-dual Riemannian
 $4$-manifolds with $H^2(Z_i, {\mathcal O}(TZ_i))=0$ for their twistor spaces.
Moreover suppose that they have positive scalar curvature.

Then,  for all sufficiently small ${\mathfrak t}>0$, the self-dual
conformal class $[g_{\mathfrak t}]$ obtained on $M_1 \# M_2$ by the
Donaldson-Friedman Theorem (\ref{thmb}) contains a metric of
positive scalar curvature.}\\

We work on the self-dual conformal classes constructed by the
Donaldson-Friedman Theorem (\ref{thmb}). Conformal Green's
Functions\cite{lom} are used to detect the sign of the scalar
curvature of these metrics. Positivity for the scalar curvature is
characterized by non-triviality of the Green's Functions. Then the
Vanishing Theorem (\ref{vanishing}) will provide the
Serre-Horrocks\cite{serremod,horrocks} vector bundle construction,
which gives the Serre Class, a substitute for the Green's Function
by Atiyah\cite{atgrn}. And non-triviality of the Serre Class will
provide the non-triviality of the extension described by it.

In sections \S\ref{construct}-\S\ref{natural} we review the
background material. In \S\ref{ssecvanishing} the vanishing theorem
is proven, and finally in \S\ref{greatsmall}-\S\ref{scalar} the sign
of the scalar curvature is detected.\\ {\bf Acknowledgments.} I want
to thank my Ph.D. advisor Claude LeBrun for his excellent
directions, letting me use his results, ideas and the figure, Justin
Sawon for his generous knowledge and many thanks to Ioana Suvaina.

\section{Self-Dual Manifolds and the Donaldson-Friedman Construction}\label{construct}

One of the main improvements in the field of self-dual Riemannian 4-manifolds is the
connected sum theorem of Donaldson and Friedman \cite{DF} published in $1989$. If $M_1$ and $M_2$ admit
self-dual metrics, then under certain circumstances their connected sum admits, too .
This helped us to create many examples of self-dual manifolds. If we state it more precisely :

\begin{thm}[Donaldson-Friedman\cite{DF}] \label{thmb}
 Let $(M_1,g_1)$ and $(M_2,g_2)$ be compact self-dual Riemannian
 $4$-manifolds and $Z_i$ denote the corresponding twistor spaces.
 Suppose that $H^2(Z_i, {\mathcal O}(TZ_i))=0$ ~~for~~ $i=1,2$.

 Then, there are self-dual conformal classes on
 $M_1 \# M_2$ whose twistor spaces arise as fibers in a
 $1$-parameter standard deformation of $Z_0=\wt{Z}_1 \cup_Q
 \wt{Z}_2.$
\end{thm}

We devote the rest of this section to understand the statement and the ideas in the
proof of this theorem since our main result (\ref{pos}) is going to be a generalization of this
celebrated theorem.

The idea is to work upstairs in the complex category rather than
downstairs. So let $p_i \in M_i$ be arbitrary points in the
manifolds. Consider their inverse images $C_i\approx \mbb{CP}_1$
under the twistor fibration, which are twistor lines, i.e. rational
curves invariant under the involution. Blow up the twistor spaces
$Z_i$ along these rational curves. Denote the exceptional divisors
by $Q_i\approx \mbb{CP}_1\times\mbb{CP}_1$ and the blown up twistor
spaces by $\wt Z_i=Bl(Z_i , C_i)$ . The normal bundles for the
exceptional divisors is computed by :

\beg{lem}[Normal Bundle]\label{normalbundle} The normal bundle of $Q_2$ in $\wt Z_2$ is
computed to be $$NQ_2 = N_{Q_2 / \wt Z_2}\approx \mco(1,-1)$$
where the second component is the fiber direction in the blowing
up process. 
\end{lem}
\beg{proof} We split the computation into the following steps
\beg{enumerate}
    \item \label{wedge} We know that $N_{C_2/Z_2}\approx\mco(1)\oplus\mco(1)$ and we compute
    its second wedge power as $$c_1(\land^2\mco(1)\oplus\mco(1))[\mbb{P}_1]=c_1(\mco(1)\oplus\mco(1))[\mbb P_1]
    =(c_1\mco(1)+c_1\mco(1))[\mbb P_1]=2$$
    by the Whitney product identity of the characteristic classes. so we have
    $$\land^2N_{C_2/Z_2}\approx\mco_{\mbb P_1}(2)$$
   \item $K_Q=\pi_1^*K_{\mbbp_1}\otimes\pi_2^*K_{\mbbp_1}=\pi_1^*\mco(-2)\otimes\pi_2^*\mco(-2)=
   \mco_{\mbbp_1\times\mbbp_1}(-2,-2)$\item  $K_Q=K_{\wt Z_2}+Q|_Q=\pi^*K_{Z_2}+2Q|_Q=
   \pi^*(K_{Z_2}|_{\mbb P_1})+2Q|_Q=
    \pi^*(K_{\mbb P_1}\otimes \land^2N_{\mbb P_1/Z_2}^*)+2Q|_Q
    =\pi^*(\mco(-2)\otimes\mco(-2))+2Q|_Q
    =\pi^*\mco(-4)+2Q|_Q$\\
since the second component is the fiber direction, the pullback
bundle will be trivial on  that so $\pi^*\mco(-4)=\mco(-4,0)$
solving for $Q|_Q$ now gives us
$$N_{Q/\wt Z_2}=Q|_Q=(K_Q\otimes\pi^*\mco(-4)^*)^{1/2}=
(\mco(-2,-2)\otimes\mco(4,0))^{1/2}=\mco(1,-1)$$

\end{enumerate}
\end{proof}

We then construct the complex analytic space $Z_0$ by identifying
$Q_1$ and $Q_2$ so that it has a normal crossing singularity
$$Z_0=\wt Z_1 \cup_Q \wt Z_2.$$
Carrying out this identification needs a little bit of care. We
interchange the components of $\bcp_1\times\bcp_1$ in the gluing
process so that the normal bundles $N_{Q_1/ \wt Z_1}$ and $N_{Q_2 /
\wt Z_2}$ are dual to each other. Moreover we should respect to the
real structures. The real structures  $\sigma_1$ and $\sigma_2$ must
agree on $Q$ obtained by identifying $Q_1$ with $Q_2$, so that the
real structures extend over $Z_0$  and form the anti-holomorphic
involution $\sigma_0 : Z_0 \to Z_0$.

Now we will be trying to deform the singular space $Z_0$, for which
the Kodaira-Spencer's standard deformation theory does not work
since it is only for manifolds it does not tell anything about the
deformations of the singular spaces. We must use the theory of
deformations of a compact reduced complex analytic spaces, which is
provided by R.\cite{friedman}. This generalized theory is quite
parallel to the theory of manifolds. The basic modification is that
the roles of $H^i(\Theta)$ are now taken up by the groups
$T^i=\tn{Ext}^i(\Omega^1, {\mathcal O})$.

We have assumed $H^2(Z_i, {\mathcal O}(TZ_i))=0$ so that the
deformations of $Z_i$ are unobstructed. Donaldson and Friedman are
able to show that $T^2_{Z_0}=\tn{Ext}^2_{Z_0}(\Omega^1, {\mathcal
O})=0$ so the deformations of the singular space is unobstructed. We
have a versal family of deformations of $Z_0$. This family is
parameterized by a neighborhood of the the origin in
$\tn{Ext}^1_{Z_0}(\Omega^1, {\mathcal O})$. The generic fiber is
non-singular and the real structure $\sigma_0$ extends to the total
space of this family.
\bct \begin{tabular}{ccc}$\omega$ : &$\mc Z \to \mc U$& ~~~for\hspace{.5cm} $Z_0=\wt Z_1 \cup_Q \wt Z_2$\\
 &$Z_0 \longmapsto 0$ &\end{tabular}\ect
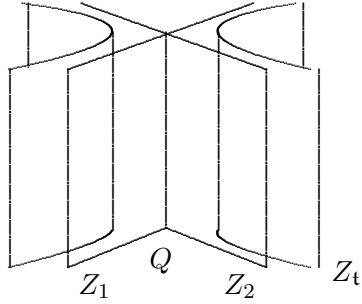
\begin{center}
\mbox{\beginpicture \setplotarea x from 0 to 200, y from -20 to 110
\put {$Z_2$} [B1] at 125  -10 \put {$Z_1$} [B1] at  70  -10 \put
{$Z_{\mathfrak t}$} [B1] at 165  -5 \put {$Q$} [B1] at    95 0
{\setlinear \plot 65 100 135 75 / \plot 60 75  130 100 / \plot 135
75 135 0 / \plot 60 75 60 0 / \plot 60 0 97 15 / \plot 97 15 135 0 /
\plot 97 89 97 15 / \plot 117 89 117  14 / \plot 155 75 155 0 /
\plot 149 100 149 76 / \plot 38 75 38 0 / \plot 77 89 77 14 / \plot
45 100 45 76 / } {\setquadratic \plot  148  75   113  89   143  100
/ \plot  148  0    120 8 113  14 / \plot 34 75 73 89 39 100 / \plot
34 0     62 8         73 14  / }\endpicture }\end{center}

Instead of working with the entire versal family, it is convenient
to work with certain subfamilies, called {\em standard
deformations}:

\beg{defn}[\cite{lom}]\index{standard deformation}  A {\em
$1$-parameter standard deformation} of $Z_0$ is a flat proper
holomorphic map $\omega : {\mathcal Z}\to {\mathcal U}\subset
\mathbb C$ of a complex $4$-manifold to an open neighborhood of
~$0$, together with an anti-holomorphic involution $\sigma :
{\mathcal Z}\to {\mathcal Z}$, such that
\begin{itemize}
\item  $\omega^{-1}(0)=Z_0$
\item $\sigma |_{Z_0}= \sigma_0$
\item $\sigma$ descents to the complex conjugation in $\mathcal U$
\item $\omega$ is a submersion away from $Q\subset Z_0$
\item $\omega$ is modeled by $(x,y,z,w)\mapsto
xy$ near any point of $Q$.\\ \end{itemize}\end{defn}We also define
\beg{defn}[Flat Map\cite{hart}]\index{flat family, map} Let $K$ be
module over a ring $A$. We say that $K$ is \emph{flat} over $A$ if
the functor $L\mapsto K \otimes_A L$ is an exact functor for all
modules $L$ over $A$.\\ Let $f:X\to Y$ be a morphism of schemes and
$\mc F$ be an $\mco_X$-module. We say $\mc F$ is \emph{flat} over
$Y$ if the stalk $\mc F_x$ is a flat $\mco_{y,Y}$-module for any
$x$. Where $y=f(x)$, $\mc F_x$ is considered to be a
$\mco_{y,Y}$-module via the natural map $f^\# : \mco_{y,Y} \to
\mco_{x,X}$. We say $X$ is \emph{flat} over $Y$ if $\mco_X$ is.
\end{defn}

Then for sufficiently small, nonzero, real ${\mathfrak t}\in
{\mathcal U}$ the complex space $Z_{\mathfrak t}=
\omega^{-1}({\mathfrak t})$ is smooth and one can show that it is
the twistor space of a self-dual metric on $M_1 \# M_2$.

\section{The Leray Spectral Sequence}\label{secleray}
Given a continuous map $f:X \to Y$ between topological spaces, and a
sheaf $\mc{F}$ over $X$, the {\em $q$-th direct image sheaf } is the
sheaf $\mathit{R}^q(f_*\mc{F})$ on $Y$ associated to the presheaf $V
\to H^q(f^{-1}(V),\mc{F}|_{f^{-1}(V)})$. This is actually the right
derived functor of the functor $f_*$. The Leray Spectral Sequence
\index{Leray Spectral Sequence} is a spectral sequence $\{E_r \}$
with
$$E^{p,q}_2 = H^p(Y,\mathit{R}^q(f_*\mc{F}))$$
$$E^{p,q}_\infty = H^{p+q}(X,\mc{F}) $$
The first page of this spectral sequence reads :
\begin{center}
\begin{tabular}{cc|cccc}
&&$\vdots$ & $\vdots$ & $\vdots$ &    \\
&&
$H^0(Y,\mathit{R}^2(f_*\mc{F}))$ &
$H^1(Y,\mathit{R}^2(f_*\mc{F}))$ &
$H^2(Y,\mathit{R}^2(f_*\mc{F}))$ &$\cdots$\\
&&
$H^0(Y,\mathit{R}^1(f_*\mc{F}))$ &
$H^1(Y,\mathit{R}^1(f_*\mc{F}))$ &
$H^2(Y,\mathit{R}^1(f_*\mc{F}))$ &$\cdots$\\
$E_2$&&
$H^0(Y,\mathit{R}^0(f_*\mc{F}))$ &
$H^1(Y,\mathit{R}^0(f_*\mc{F}))$ &
$H^2(Y,\mathit{R}^0(f_*\mc{F}))$ &$\cdots$\\
\cline{3-6}
\end{tabular}
\end{center}\vspace{.3cm}
A degenerate case is when $\mathit{R}^i(f_*\mc{F})=0 $ for all $i
> 0$.
\begin{rmk} This is the case if $\mc{F}$ is flabby for
example. Remember that to be {\em flabby} means that the restriction
map $r: \mc{F}(\mathit{B}) \to \mc{F}(\mathit{A})$ is onto for open
sets $\mathit{B} \subset \mathit{A}$. In this case $H^i(X,\mc{F})=0$
for $i > 0$ as well as $H^i(\mathit{U},\mc{F}|_\mathit{U})=0$ for
$\mathit{U}$ open, because the restriction of a flabby sheaf to any
open subset is again flabby by definition. That means
$H^q(f^{-1}(.),\mc{F}|_.)=0$ for all $q > 0$ so
$\mathit{R}^i(f_*\mc{F})=0 $. \end{rmk}

When the spectral sequence degenerates this way, the second and
succeeding rows of the first page vanish. And because $V \to
H^0(f^{-1}(V),\mc{F}|_{f^{-1}(V)})$ is the presheaf of the direct
image sheaf, we have $\mathit{R}^0f_*=f_*$. So the first row consist
of $H^i(Y,f_*\mc{F})$'s. Vanishing of the differentials cause
immediate convergence to $E^{i,0}_\infty=H^{i+0}(X,\mc{F})$. So we
got:
\begin{prop} \label{leray} If $\mathit{R}^i(f_*\mc{F})=0 $ for all $i>0$, then
$H^i(X,\mc{F})=H^i(Y,f_*\mc{F})$ naturally for all $i \geq 0$.
\end{prop}

As another example, the following proposition reveals a different
sufficient condition for this degeneration. See \cite{voisin} v2 ,
p124 for a sketch of the proof:

\begin{prop}[Small Fiber Theorem] \label{smallfiber} Let $f:X \to Y$ be a holomorphic, proper and
submersive map between complex manifolds, $\mc{F}$ a coherent
analytic sheaf or a holomorphic vector bundle on $X$. Then
$H^i(f^{-1}(y),\mc{F}|_{f^{-1}(y)})=0$ for all $y \in Y$ implies
that $\mathit{R}^i(f_*\mc{F})=0$.
\end{prop}

As an application of these two propositions, we obtain the main
result of this section :

\begin{prop} \label{lerayblowup} Let $Z$ be a complex $n$-manifold with a complex
$k$-dimensional submanifold $V$. Let $\wt{Z}$ denote the blow up of
$Z$ along $V$, with blow up map $\pi : \wt{Z} \to Z$. Let $\mc{G}$
denote a coherent analytic sheaf(or a vector bundle) over $Z$. Then
we can compute the cohomology of $\mc{G}$ on either side i.e.
$$H^i(\wt{Z},\pi^*\mc{G})=H^i(Z,\mc{G}).$$
\end{prop}
\begin{proof} The inverse image of a generic point on $Z$ is a point,
else a $\mathbb{P}^{n-k-1}$. We have
$$H^i(f^{-1}(y),\pi^*\mc{G}|f^{-1}(y))=H^i(\mathbb{P}^{n-k-1},\mc{O})=H^{0,i}_{
\bar{\partial }}(\mathbb{P}^{n-k-1})=0$$ at most, since the
cohomology of $\mathbb{P}^{n-k-1}$ is accumulated in the middle for
$i>0$. So that we can apply Proposition (\ref{smallfiber}) to get
$\mathit{R}^i(\pi_*\pi^\ast \mc{G})=0 $ for all $i>0$. Which is the
hypothesis of Proposition (\ref{leray}), so we get
$H^i(\wt{Z},\pi^*\mc{G})=H^i(Z,\pi_*\pi^*\mc{G}) $ naturally for all
$i\geq0$, and the latter equals $H^i(Z,\mc{G})$ since
$\pi_*\pi^*\mc{G}=\mc{G}$ by the combination of the following two
lemmas.\end{proof}

\beg{lem}[Projection Formula\cite{hart}] If $f: (X,\mc{O}_X) \to
(Y,\mc{O}_Y)$ is a morphism of ringed spaces, if $\mc{F}$ is an
$\mc{O}$-module, and if $\mc{E}$ is a locally free $\mc{O}_Y$-module
of finite rank. Then there is a natural isomorphism
$$f_*(\mc{F}\otimes_{\mc{O}_X}f^*\mc{E})=f_*\mc{F}\otimes_{\mc{O}_Y}\mc{E},$$
in particular for $\mc{F}=\mc{O}_X$
$$f_*f^*\mc{E}=f_*\mc{O}_X\otimes_{\mc{O}_Y}\mc{E}.$$
\end{lem}

\begin{lem}[Zariski's Main Theorem,Weak Version\cite{hart}]
Let $f:X \to Y$ be a birational projective morphism of noetherian
integral schemes, and assume that $Y$ is normal. Then
$f_*\mc{O}_X=\mc{O}_Y$.
\end{lem}

\begin{proof} The question is local on $Y$. So we may assume that
$Y$is affine and equal to $SpecA$. Then $f_*\mc{O}_X$ is a coherent
sheaf of $\mc O_Y$-algebras, so $B=\Gamma(Y,f_*\mc{O}_X)$ is a
finitely generated $A$-module. But $A$ and $B$ are integral domains
with the same quotient field, and $A$ is integrally closed, we must
have $A=B$. Thus $f_*\mc{O}_X=\mc{O}_Y$
\end{proof}


\section{Natural square root of the canonical bundle of a twistor
space}\label{natural}

In the next section, we are going to prove that a certain cohomology
group of a line bundle vanishes. For that we need some lemmas. First
of all, the canonical bundle of a twistor space $Z$ has a natural
square root, equivalently $Z$ is a spin manifold as follows:

The Riemannian connection of $M$ acts on the $2$-forms hence the
twistor space, accordingly we can split the tangent bundle
$T_xZ=T_xF \oplus (p^*TM)_x$. The complex structure on $(p^*TM)_x$
is obtained from the identification $\cdot\varphi : T_xM
\longleftrightarrow (\mbb{V}_+)_x$ provided by the Clifford
multiplication of a non-zero spinor $\varphi \in (\mbb{V}_+)_x$.
This identification is linear in $\varphi$ as $\varphi$ varies over
$({\mbb{V}_+})_x$. So we have a nonvanishing section of
$\mco_Z(1)=\mco_{\mbb P \mbb V_-}(1)$ with values in
$Hom(TM,{\mathbb V}_+)$ or $Hom(p^*TM,p^*{\mbb V}_+)$ trivializing
the bundle $$\mco_Z(1)\otimes Hom(p^*TM,p^*{\mbb V}_+)
=\mco_Z(1)\otimes {p^*TM}^* \otimes p^*{\mbb V}_+ \approx {p^*TM}^*
\otimes \mco_Z(1)\otimes p^*{\mbb V}_+$$ hence yielding a natural
isomorphism
\begin{equation}{p^*TM}\approx \mco_Z(1)\otimes p^*{\mbb V}_+ ,\label{iso1}
\end{equation}
where $\mco_Z(1)=\mco_{\mbb P \mbb V_-}(1)$ is the positive Hopf
bundle on the fiber.

The Hopf bundle exist locally in general, so as the isomorphism. If
$M$ is a spin manifold $\mbb V_\pm$ exist globally on $M$ and
$\mco_Z(1)$ exist on $Z$, so our isomorphism holds globally.
Furthermore, we have a second isomorphism holding for any projective
bundle, obtained as follows (see \cite{fultoni} p434 , \cite{zheng}
p108) :

Let E be a complex vector bundle of rank-$(n+1)$ over $M$, $p: \mbb
PE \to M$ its projectivization. We have the imbedding of the
tautological line bundle $\mco_{\mbb PE}(-1)\hookrightarrow p^*E$.
Giving the exact sequence
$$0 \to \mco_{\mbb PE}(-1) \to p^*E \to{ {p^*E}/{\mco_{\mbb PE}(-1)} }\to 0,$$
tensoring by $\mco_{\mbb PE}(1)$ gives
$$0 \to \mco_{\mbb PE} \to \mco_{\mbb PE}(1)\otimes p^*E \to T_{\mbb PE/M} \to 0$$
where $T_{\mbb PE/M} \approx
Hom(\mco(-1),\mco(-1)^\perp)=Hom(\mco(-1),p^*E/{\mco(-1)})=\mco(1)\otimes
{p^*E/\mco(-1)} $  is the relative tangent bundle of $\mbb P E$
over $M$, originally defined to be ${\Omega^1_{\mbb PE/M}}^*$.
Taking $E={\mbb V}_-$ , $TF$ denoting the tangent bundle over the
fibers:
$$0 \to \mco_Z \to \mco_Z(1)\otimes p^*\mbb V_- \to TF \to 0$$
so, we got our second isomorphism :
\begin{equation} TF\oplus\mco_Z \approx \mco_Z(1)\otimes
p^*\mbb V_- \label{iso2} \end{equation}

Now we are going to compute the first chern class of the spin
bundles $\mbb V_\pm$, and see that $c_1(\mbb V_\pm)=0$. Choose a
connection $\nabla$ on $\mbb V_\pm$. Following \cite{kn} it defines
a connection on the associated principal $\mathfrak{su}(2)$ bundle
$P$, with connection one form $\omega \in A^1(P,\mf{su}(2))$ defined
by the projection \cite{morita} $T_uP \to V_u \approx \mf{su}(2)$
having curvature two form $\Omega \in A^2(P,\mf{su}(2))$ defined by
\cite{kn} :
$$\Omega(X,Y)=d\omega(X,Y)+\frac{1}{2}[\omega(X),\omega(Y)]\ \mathrm{ for } \ X,Y \in
T_uP.$$ We define the first polynomial functions $f_0,f_1,f_2$ on
the lie algebra $\mf{su}(2)$ by
$$det(\lambda I_2+\frac{i}{2\pi}M)=\sum_{k=0}^2{f_{2-k}(M)\lambda^k}=f_0(M)\lambda^2
+f_1(M)\lambda+f_2(M) \ \mathrm{for} \ M \in \mf{su}(2).$$ Then
these polynomials $f_i: \mf{su}(2)\to \mbb R$ are invariant under
the adjoint action of $SU(2)$, denoted $f_i \in I^1(SU(2))$, namely
$$f_i(\ad_g(M))=f_i(M) \ \mathrm{for} \ g \in SU(2) \ , M \in
\mf{su}(2)$$ where $ \ad_g: \mf{su}(2) \to \mf{su}(2)$~ is defined
by
~$\ad_g(M)=R_{g^{-1}*}L_{g*}(M)$\label{adjointrepresentation}\index{adjoint
representation}.\\ If we apply any $f \in I^1(SU(2))$ after $\Omega$
we obtain:
$$f \circ \Omega \ : \ T_uP \times T_uP \lra \mf{su}(2) \lra \mbb{R}.$$
It turns out that $f \circ \Omega$ is a closed form and projects
to a unique $2$-form say $\overline{f \circ \Omega}$ on $M$ i.e.
$f \circ \Omega=\pi^*(\overline{f \circ \Omega})$ where $\pi:P\to
M$. By the way, a $q$-form $\varphi$ on $P$ projects to a unique
$q$-form, say $\overline{\varphi}$ on $M$ if $\varphi(X_1
\cdot\cdot \, X_q)=0$ whenever at least one of the $X_i$'s is
vertical and $\varphi(R_{g*}X_1 \cdot\cdot \,
R_{g*}X_q)=\varphi(X_1 \cdot\cdot \, X_q)$. $\overline{\varphi}$
on $M$ defined by $\overline{\varphi}(V_1 \cdot\cdot \,
Vq)=\varphi(X_1 \cdot\cdot \, X_q) , \pi(X_i)=V_i$ is independent
of the choice of $X_i$'s. See \cite{kn}v2p294 for details.\\
So, composing with $\Omega$ and projecting defines a map w
$:I^1(SU(2)) \to H^2(M,\mbb R)$ called the Weil homomorphism, it
is actually an algebra homomorphism when extended to the other
gradings.\\
Finally, the chern classes are defined by $c_k(\mbb
V_\pm):=\left[\overline{f_k \circ \Omega}\right]$ independent of
the connection chosen. Notice that
$f_2(M)=det(\frac{i}{2\pi}M),f_1(M)=tr(\frac{i}{2\pi}M)$ in our
case. And if $M\in \mf{su}(2)$ then $e^M \in SU(2)$ implying
$1=det(e^M)=e^{trM}$ and $trM=0$. But $\Omega$ is of valued
$\mf{su}(2)$, so if you apply the $f_1=tr$ after $\Omega$ you get
$0$. Causing $c_1(\mbb V_\pm)=\left[\overline
{tr(\frac{i}{2\pi}\Omega)}\right]=0$. \\
One last remark is that $\overline{f_k \circ \Omega}=\gamma_k$ in
the notation of \cite{kn},
$\gamma_1=P^1(\frac{i}{2\pi}\Theta)=tr(\frac{i}{2\pi}\Theta)$ in
the notation of \cite{gh}p141,p407. And $\Omega=\pi^*\Theta$ in
the line bundle case.

Vanishing of the first chern classes mean that the determinant line
bundles of $\mbb V_\pm$ are diffeomorphically trivial since
$c_1(\wedge^2\mbb V_\pm)=c_1\mbb V_\pm=0$. Combining this with the
isomorphisms (\ref{iso1}) and (\ref{iso2}) yields: \bct
$\wedge^2p^*TM=\wedge^2(\mco_Z(1)\otimes p^*\mbb V_+)=
\mco_Z(2)\otimes\wedge^2p^*\mbb
V_+=\mco_Z(2)=\mco_Z(2)\otimes\wedge^2p^*\mbb V_- $
$=\wedge^2(\mco_Z(1)\otimes p^*\mbb V_-)=\wedge^2(TF \oplus
\mco_Z)=\oplus_{2=p+q}(\wedge^pTF\otimes\wedge^q\mco_Z)=TF\otimes\mco_Z=TF$
\ect
since $TF$ is a line bundle. Taking the first chern class of both
sides: $$ c_1(p^*TM)=c_1(\wedge^2p^*TM)=c_1TF. \label{3} $$
Alternatively, this chern class argument could be replaced with the
previous taking wedge powers steps if the reader feels more
comfortable with it.\\
Last equality implies the decomposition:
$$c_1Z=c_1(p^*TM\oplus TF)=c_1(p^*TM)+c_1TF=2c_1TF.$$ So, $TF^*$ is a
differentiable square root for the canonical bundle of $Z$. If $M$
is not spin $\mbb V_\pm,\mco_Z(1)$ are not globally defined, but
the complex structure on their tensor product is still defined,
and we can still use the isomorphisms (\ref{iso1}),(\ref{iso2})
for computing chern classes of the almost complex structure on $Z$
using differential forms defined locally by the metric.
Consequently our decomposition is valid whether $M$ is spin or
not.

One more word about the differentiable square roots is in order here.
A differentiable square root implies a holomorphic one on complex
manifolds since in the sheaf sequence:
\begin{center}
\begin{tabular}{ccccccc}
.. & $\to$ & $H^1(M,\mco^*)$ & $\to$ & $H^2(M,\mbb Z)$ & $\to$ & $..$ \\
   &       &    $L$          & $\mapsto$ & $c_1(L)$ &  $\mapsto$ & $0$ \\
   &       &           &  & $\frac{1}{2}c_1(L)$ &  $ \mapsto$ & $0$  \end{tabular}
\end{center}
$c_1(L)$ maps to $0$ since it is coming from a line bundle, and if it
decomposes, $\frac{1}{2}c_1(L)$ maps onto $0$ too, that means it is
the first chern class of a line bundle.

\section{Vanishing Theorem}\label{ssecvanishing}

Let $\omega : \mathcal Z \to \mathcal U$ be a $1$-paramater standard deformation of $Z_0$, where
$\mc U \subset \mbb C$ is an open disk about the origin. Then the invertible sheaf $K_{\mc Z}$
has a square root as a holomorphic line bundle as follows:

We are going to show that the Steifel-Whitney class $w_2(K_{\mc Z})$
is going to vanish. We write $\mc Z={\mc U}_1 \cup {\mc U}_2$ where
${\mc U}_i$ is a tubular neighborhood of $\wt{Z}_i$, ${\mc U}_1 \cap
{\mc U}_2$ is a tubular neighborhood of $Q=\wt{Z}_1 \cap \wt{Z}_2$.
So that ${\mc U}_1 , {\mc U}_2$ and ${\mc U}_1 \cap {\mc U}_2$
deformation retracts on $\wt{Z}_1 , \wt{Z}_2$ and $Q$. Since $Q
\approx {\mbb P}_1 \times {\mbb P}_1$ is simply connected, $H^1({\mc
U}_1 \cap {\mc U}_2,{\mbb Z}_2)=0$ and the map $r_{12}$ in the
Mayer-Vietoris exact sequence : \bct
\begin{tabular}{ccccccc}
\llap{ .. ~~$\to$~~~}$H^1({\mc U}_1 \cap {\mc U}_2,{\mbb Z}_2)$ &
$\to$ & $H^2({\mc U}_1 \cup {\mc U}_2,{\mbb Z}_2)$ &
$\stackrel{r_{12}}{\to}$ &
$H^2({\mc U}_1,{\mbb Z}_2) \oplus H^2({\mc U}_2,{\mbb Z}_2)$ 
& $\to$ &.. \\
$\parallel$ && $\in$ &  &&& \\
$0$  && $w_2(K_{\mc Z})$ &&&& \end{tabular} \ect is injective.
Therefore it is enough to see that the restrictions $r_i(w_2(K_{\mc
Z})) \in H^2({\mc U}_i,{\mbb Z}_2)$ are zero. For that, we need to
see that $K_{\mc Z}|_{\wt{Z}_i}$ has a radical : \bct
$K_\mcz|_{\wt{Z}_1} \stackrel{(1)}{=}
(K_{\wt{Z}_1}-\wt{Z}_1)|_{\wt{Z}_1} \stackrel{(2)}{=}
(K_{\wt{Z}_1}+Q)|_{\wt{Z}_1} \stackrel{(3)}{=}
((\pi^*K_{Z_1}+Q)+Q)|_{\wt{Z}_1}=
2(\pi^*K^{1/2}_{Z_1}+Q)|_{\wt{Z}_1}$ \ect where $(1)$ is the
application of the adjuction formula on $\wt{Z}_1$,
$ K_{\wt{Z}_1} = K_{\mcz}|_{\wt{Z}_1} \otimes [\wt{Z}_1] $.\\
$(2)$ comes from the linear equivalence of $0$ with $Z_t$ on $\wt{Z}_1$,
and $Z_t$ with $Z_0$:
$$0=\mco(Z_t)|_{\wt{Z}_1}=\mco(Z_0)|_{\wt{Z}_1}=\mco(\wt{Z}_1+\wt{Z}_2)|_{\wt{Z}_1}
=\mco(\wt{Z}_1+Q)|_{\wt{Z}_1}$$
(3) is the change of the canonical bundle under the blow up along a submanifold, see \cite{gh}p608.
$K_{Z_1}$ has a natural square root as we computed in the previous section, so
$\pi^*K^{1/2}_{Z_1}\otimes[Q]$ is a square root of $K_\mcz$ on $\widetilde{Z}_1$. Similarly
on $\widetilde{Z}_2$, so $K_\mcz$ has a square root $K^{1/2}_\mcz$.

Before our vanishing theorem, we are going to mention the Semicontinuity
Principle and the
Hitchin's Vanishing theorem, which are involved in the proof:

\begin{lem}[Semicontinuity Principle\cite{voisin}]
\label{semicontinuity}
Let $\phi : \mc X \to \mc B$ be a family of complex compact manifolds
With fiber $X_b , b \in \mc B$. Let $\mc F$ be a holomorphic vector bundle over
$\mc X$, then

The function $b \mapsto h^q(X_b,\mc F|_{X_b})$ is upper semicontinuous.
In other words, we have $ h^q(X_b,\mc F|_{X_b}) \leq h^q(X_0,\mc F|_{X_0})$ for $b$ in a
neighborhood of $0 \in \mc B$.
\end{lem}

\begin{lem}[Hitchin Vanishing\cite{hitlin}\cite{poon1}]
\label{hitchinvanishing}
Let $Z$ be the twistor space of an oriented self-dual riemannian manifold
of positive scalar curvature with canonical bundle $K$, then
$$ h^0(Z,\mco(K^{n/2}))=h^1(Z,\mco(K^{n/2}))=0 ~ for ~ all ~ n \geq 1.$$

\end{lem}

\begin{thm}[Vanishing Theorem]
\label{vanishing} \index{Vanishing Theorem}
Let $\omega : {\mathcal Z} \to {\mathcal U}$
be a $1$-parameter standard deformation of $Z_0$, where $Z_0$ is as
in Theorem (\ref{thmb}), and ${\mathcal U}\subset \mathbb C$ is a
neighborhood of the origin. Let $L\to {\mathcal Z}$ be the
holomorphic line bundle defined by
$${\mathcal O}(L^*)
={\mathscr I}_{\wt{Z}_1}(K_{\mathcal Z}^{1/2})$$
If $(M_i,[g_i])$ has positive scalar curvature,
then  by possibly replacing  ${\mathcal U}$ with a smaller
 neighborhood  of $0\in \mathbb C$ and simultaneously replacing
 ${\mathcal Z} $ with its inverse image, we can arrange for
 our complex $4$-fold $\mathcal Z$ to satisfy
 $$H^1 ({\mathcal Z} , {\mathcal O}(L^*))=H^2 ({\mathcal Z} , {\mathcal
   O}(L^*))=0.$$
\end{thm}

\begin{proof} The proof proceeds by analogy to the techniques in \cite{lom}, and consists of several steps :\beg{enumerate}

\item It is enough to \textbf{show that $\mathbf{H^j(Z_0,\mco(L^*))=0}$ for $\mathbf{j=1,2}$ 
:} \\ Since that would imply $h^j(Z_t,\mco(L^*)) \leq 0$ for $j=1,2$
in a neighborhood by the semicontinuity principle. Intuitively, this
means that the fibers are too small, so we can apply Proposition
(\ref{smallfiber}) to see $R^j\omega_*\mco(L^*)=0$ for $j=1,2$. The
first page of the Leray Spectral Sequence reads :\newline
\begin{center}
\begin{tabular}{|ccccc}
 $\vdots$ & $\vdots$ & $\vdots$ &    \\
$H^0(\mc U,R^3\omega_*\mco(L^*))$ &
$H^1(\mc U,R^3\omega_*\mco(L^*))$ &
$H^2(\mc U,R^3\omega_*\mco(L^*))$ &$\cdots$\\
$0$ &
$0$ &
$0$ &$\cdots$\\
$0$ &
$0$ &
$0$ &$\cdots$\\

\llap{$E_2$\hspace{.5cm} }$H^0(\mc U,R^0\omega_*\mco(L^*))$ &
$H^1(\mc U,R^0\omega_*\mco(L^*))$ &
$H^2(\mc U,R^0\omega_*\mco(L^*))$ &$\cdots$\\
\cline{0-3}
\end{tabular}
\end{center}\vspace{.9cm}
Remember that $$E^{p,q}_2 = H^p(\mc U,R^q\omega_*\mco(L^*))$$
$$E^{p,q}_\infty = H^{p+q}(\mc Z,\mco(L^*)) $$ and
that the differential $$d_2(E^{p,q}_2) \subset E^{p+2,q-1}_2.$$
Vanishing of the second row implies the immediate convergence of the
first row till the third column because of the differentials, so

$$E^{p,0}_\infty=E^{p,0}_2 ~\textnormal{i.e.}~ H^{p+0}(\mc
Z,\mco(L^*))=H^p(\mc U,R^0\omega_*\mco(L^*)) ~\textnormal{for}~ p\leq3$$
hence $H^p(\mc Z,\mco(L^*))=H^p(\mc U,R^0\omega_*\mco(L^*))$ , for
$p\leq3$.\\
Since $\mc U$ is one dimensional, $\omega : \mc Z \to \mc U$ has to
be a flat morphism, so the sheaf $\omega_*\mco(L^*)$ is
coherent\cite{gunningr3,bonica}. $\mc U$ is an open subset of $\mbb
C$ implying that it is Stein. And the so called Theorem B of Stein
Manifold theory characterizes them as possesing a vanishing higher
dimensional($p>0$) coherent sheaf cohomology \cite{lewis}p67,
\cite{hart}p252, \cite{gunningr3,bonica}. So $H^p(\mc
U,\omega_*\mco(L^*))=0$ for $p>0$. Tells us that $H^p(\mc
Z,\mco(L^*))=0$ for $0<p\leq3$.

\item \label{mayervietoris}Related to $Z_0$, we have the \textbf{Mayer-Vietoris like} sheaf exact
sequence
$$0\to {\mathcal O}_{Z_0}(L^*)\to \nu_*{\mathcal O}_{\wt{Z}_1}(L^*)\oplus
\nu_*{\mathcal O}_{\wt{Z}_2}(L^*) \to {\mathcal O}_Q(L^*)\to 0$$
where $\nu : \wt{Z}_1\sqcup \wt{Z}_2\to Z_0$ is the inclusion map on
each of the two components of the disjoint union $\wt{Z}_1\sqcup
\wt{Z}_2$. This gives the long exact cohomology sequence piece :
\bct $0\to H^1(\mco_{Z_0}(L^*))\to H^1(Z_0,\nu_*\mco_{\wt{Z}_1}(L^*)
\oplus \nu_*\mco_{\wt{Z}_2}(L^*)) \to H^1(\mco_Q(L^*))\to
H^2(\mco_{Z_0}(L^*))\to H^2(Z_0,\nu_*\mco_{\wt{Z}_1}(L^*)\oplus
\nu_*\mco_{\wt{Z}_2}(L^*))\to 0$ \ect due to the fact that :

\item \label{restrto quadric cohomology} $\mathbf{H^0(\mco_Q(L^*))= H^2(\mco_Q(L^*))=0}$ \textbf{:} To see this, we have
to understand the restriction of $\mco(L^*)$ to $Q$ :   \bct
$L^*|_Q=(\frac{1}{2}K_\mc{Z}-\wt{Z}_1)|_{\wt{Z}_2}|_Q =(
\frac{1}{2}( K_{\wt{Z}_2}-\wt{Z}_2 )-\wt{Z}_1  )  |_{\wt{Z}_2}|_Q =(
\frac{1}{2}( K_{\wt{Z}_2} +  Q     )-  Q       )  |_{\wt{Z}_2}|_Q
=\frac{1}{2}( K_{\wt{Z}_2} -  Q     ) |_{\wt{Z}_2}|_Q =\frac{1}{2}(
K_Q -  Q  -  Q  )  |_{\wt{Z}_2}|_Q =(\frac{1}{2} K_Q  -  Q       )
|_{\wt{Z}_2}|_Q =\frac{1}{2} K_Q |_Q \otimes NQ_{\wtz_2}^{-1}
=\mco(-2,-2)^{1/2} \otimes {\mco(1,-1)}^{-1}=\mco(-2,0)$  \ect
here,  we have computed the normal bundle of $Q$ in $\wt{Z}_2$ in
Lemma (\ref{normalbundle}) as $\mco(1,-1)$, where the second
component is the fiber direction in the blowing up process. So the
line bundle $L^*$ is trivial on the fibers. Since $Q={\mbb P}_1
\times {\mbb P}_1$, we have \bct $H^0({\mbb P}_1 \times {\mbb
P}_1,\mco(-2,0))=H^0({\mbb P}_1 \times {\mbb P}_1,\pi_1^*\mco(-2))
=H^0({\mbb P}_1 ,{\pi_1}_*\pi_1^*\mco(-2))=H^0({\mbb
P}_1,\mco(-2))=0$ \ect by the Leray spectral sequence and the
projection formula since $H^k({\mbb P}_1,\mco)=0$ for $k>0$.
Similarly \bct $H^2({\mbb P}_1 \times {\mbb
P}_1,\mco(-2,0))=H^2({\mbb P}_1,\mco(-2))=0$\ect by dimensional
reasons. Moreover, for the sake of curiosity \bct $H^1({\mbb P}_1
\times {\mbb P}_1,\mco(-2,0))=H^1({\mbb P}_1,\mco(-2)) \approx
H^0({\mbb P}_1,\mco(-2)\otimes \mco(-2)^*)^*=H^0({\mbb
P}_1,\mco)^*=\mbb C$.\ect

\item \label{hitvan}
$\mathbf{H^1(\wtz_2,\mco_{\wtz_2}(L^*))=H^2(\wtz_2,\mco_{\wtz_2}(L^*))=0}$
\textbf{:} These are applications of Hitchin's second Vanishing
Theorem and are going to help us to simplify our exact sequence
piece. \bct
$H^1(\wtz_2,\mco_{\wtz_2}(L^*))=H^1(\wtz_2,\mco(K_\mcz^{1/2}-\wtz_1)|_{\wtz_2})
=H^1(\wtz_2,\mco(K_\mcz^{1/2}-Q)|_{\wtz_2})=H^1(\wtz_2,\pi^*K_{Z_2}^{1/2})
=H^1(Z_2,\pi_*\pi^*K_{Z_2}^{1/2})=H^1(Z_2,K_{Z_2}^{1/2})=0$  \ect by
the Leray spectral sequence, projection formula and Hitchin's
Vanishing theorem for $Z_2$, since it is the twistor space of a
positive scalar curvature space. This implies
$H^2(Z_2,K_{Z_2}^{1/2}) \approx H^1(Z_2,K_{Z_2}^{1/2})^*=0$ because
of the Kodaira-Serre Duality. Hence our cohomological exact sequence
piece simplifies to \bct $0\to H^1(\mco_{Z_0}(L^*))\to
H^1(\wtz_1,\mco_{\wt{Z}_1}(L^*) ) \to H^1(\mco_Q(L^*))\to
H^2(\mco_{Z_0}(L^*))\to H^2(\wtz_1,\mco_{\wt{Z}_1}(L^*) )\to 0$ \ect

\item \label{tech}\textbf{$\mathbf{H^k(\mco_{\wtz_1}(L^* \otimes [Q]^{-1}_{\wtz_1}) )=0}$ for
$\mathbf{k=1,2,3}$ :} This technical result is going to be needed to
understand the exact sequence in the next step. First we simplify
the sheaf as \bct
$(L^*-Q)|_{\wtz_1}\stackrel{def}=(\frac{1}{2}K_\mcz-\wtz_1-Q)|_{\wtz_1}
=\frac{1}{2}K_\mcz|_{\wtz_1}\stackrel{adj}=\frac{1}{2}(K_{\wtz_1}-\wtz_1)|_{\wtz_1}
=\frac{1}{2}(K_{\wtz_1}+Q)|_{\wtz_1}.$ \ect So \bct
$H^k(\wtz_1,L^*-Q)=H^k(\wtz_1,(K_{\wtz_1}+Q)/2 )\stackrel{sd}\approx
H^{3-k}(\wtz_1,(K_{\wtz_1}-Q)/2)^*$
$=H^{3-k}(\wtz_1,\frac{1}{2}\pi^*K_{Z_1})^*
\stackrel{lss}=H^{3-k}(Z_1,\frac{1}{2}\pi_*\pi^*K_{Z_1})^*$
$\stackrel{pf}=H^{3-k}(Z_1,K_{Z_1}^{1/2})^*\stackrel{sd}\approx
H^k(Z_1,K_{Z_1}^{1/2})$\ect and one of the last two terms vanish in
any case for $k=1,2,3$. So we apply the Hitchin Vanishing theorem
for dimensions $0$ and $1$.

\item \label{restrictiontoq}\textbf{Restriction maps to $\mathbf{Q}$ :} Consider the exact sequence of
sheaves on $\wtz_1$ :
$$0 \to \mco_{\wtz_1}(L^* \otimes [Q]^{-1}_{\wtz_1}) \to
\mco_{\wtz_1}(L^*) \to \mco_Q(L^*) \to 0.$$ The previous step 
implies that the restriction maps :
$$H^1(\mco_{\wtz_1}(L^*))\stackrel{restr_1}{\longrightarrow}H^1(\mco_Q(L^*))$$
and
$$H^2(\mco_{\wtz_1}(L^*))\stackrel{restr_2}{\longrightarrow}H^2(\mco_Q(L^*))$$
are isomorphism. In particular $H^2(\mco_{\wtz_1}(L^*))=0$ due to
(\ref{restrto quadric cohomology}). 
Incidentally, this exact sheaf sequence is a substitute for the role
played by the Hitchin Vanishing Theorem, for the $\wtz_2$ components
in the cohomology sequence. It also assumes Hitchin's theorems for
the $\wtz_1$ component.

\item\textbf{Conclusion :} Our cohomology exact sequence piece reduces to
\bct $0\to H^1(\mco_{Z_0}(L^*))\to H^1(\wtz_1,\mco_{\wt{Z}_1}(L^*) )
\stackrel{restr_1}{\longrightarrow} H^1(\mco_Q(L^*))\to
H^2(\mco_{Z_0}(L^*))\to 0$ \ect the isomorphism in the middle forces
the rest of the maps to be $0$ and hence we get
$H^1(\mco_{Z_0}(L^*))=H^2(\mco_{Z_0}(L^*))=0$.
\end{enumerate}
\end{proof}

\section*{The Sign of the Scalar Curvature}\label{secsign}

The sections after this point are devoted to detect the sign of the
scalar curvature of the metric we consider on the connected sum. We
use Green's Functions for that purpose. Positivity for the scalar
curvature is going to be characterized by non-triviality 
of the Green's
Functions. Then our Vanishing Theorem will provide the
Serre-Horrocks vector bundle construction, which gives the Serre
Class, a substitute for the Green's Function by Atiyah\cite{atgrn}.
And non-triviality 
of the Serre Class will provide the non-triviality
of the extension described by it.

\section{Green's Function Characterization}\label{greatsmall}

In this section, we define the Green's Functions. To get a unique Green's Function, we need an operator
which has a trivial kernel. So we begin with a compact Riemannian $4$-manifold $(M,g)$,
 and assume that
its \emph{Yamabe Laplacian}\index{Laplacian, Yamabe} $\Delta + s/6$ has trivial kernel.This is
automatic if $g$ is conformally equivalent to a metric of positive 
scalar curvature, impossible if it is conformally equivalent to a
metric of zero scalar curvature because of the Hodge Laplacian, and
may or may not happen for a metric of negative scalar curvature.
Since the Hodge Laplacian $\Delta$ is self-adjoint, $\Delta+s/6$ is
also self-adjoint implying that it has a trivial cokernel, if once
have a trivial kernel. Therefore it is a bijection and we have a
unique smooth solution $u$ for the equation ~$(\Delta + s/6) u =f$~
for any smooth function $f$. It also follows that it has a unique
distributional solution $u$ for any distribution $f$. Let $y \in M$
be any point. Consider the Dirac delta distribution $\delta_y$ at
$y$ defined by\index{Dirac delta distribution}
$$\delta_y : C^\infty(M)\to \mbb R ~~,~~ \delta_y(f)=f(y)$$
intuitively, this behaves like a function identically zero on
$M-\{y\}$, and infinity at $y$ with integral $1$. Then there is a
unique distributional solution $G_y$ to the equation
$$(\Delta + s/6)G_y =\delta_y$$ called the \emph{Green's Function}
\index{Green's function} for $y$.
Since $\delta_y$ is identically zero on $M-\{ y\}$, elliptic
regularity implies that $G_y$ is smooth on $M-\{y\}$.

About $y$, one has an expansion
$$G_y =  \frac{1}{4\pi^2}\frac{1}{\varrho^2_y}+ O(\log \varrho_y)$$
near  $\varrho_y$ denotes the distance from $y$. In the case $(M,g)$
is self-dual this expansion reduces to \cite{atgrn}
$$G_y =\frac{1}{4\pi^2}\frac{1}{\varrho_y^2}+bounded~ terms$$
We also call $G_y$ to be the {\em conformal Green's function} of
$(M,g,y)$.

This terminology comes from the fact that the
 Yamabe Laplacian is  a {\em conformally invariant}
differential operator as a map between sections of some real line
bundles. For any nonvanishing smooth function $u$, the conformally
equivalent metric $\tilde{g}=u^2g$ has scalar curvature
$$\tilde s=6u^{-3}(\Delta+s/6)u$$
A consequence of this is that $u^{-1}G_y$ is the conformal Green's
function for $(M,u^2g,y)$ if $G_y$ is the one for $(M,g,y)$.

Any metric on a compact manifold is conformally equivalent to a
metric of constant scalar curvature sign. Since if $u\not\equiv 0$
is the eigenfunction of the lowest eigenvalue $\lambda$ of the
Yamabe Laplacian, $$\tilde{s}=6u^{-3}\lambda u=6\lambda u^{-2}$$ for
the metric $\tilde g=u^2g$. Actually a more stronger statement is
true thanks to the proof\cite{lp} of the Yamabe Conjecture, any
metric on a compact manifold is conformally equivalent to a metric
of constant scalar curvature(CSC). Also if two metrics with scalar
curvatures of fixed signs are conformally equivalent, then their
scalar curvatures have the same sign.

The sign of Yamabe constant of a conformal class, meaning the sign
of the constant scalar curvature of the metric produced by the proof
of the Yamabe conjecture is the same as the sign of the smallest
Yamabe eigenvalue $\lambda$ for any metric in the conformal class.

Before giving our characterization for positivity, we are going to
state the maximum principle we will be using. Consider the
differential operator $L_c =
\sum_{i,j=1}^na^{ij}(x)\frac{\partial^2}{\partial x_i
\partial x_j}$ arranged so that $a^{ij}=a^{ji}$. It is called
\emph{elliptic \cite{protter}}\index{ellipticity}
\label{ellipticity} at a point $x=(x_1 .. x_n)$ if there is a
positive quantity $\mu (x)$ such that
$$\sum_{i,j=1}^na^{ij}(x)\xi_i\xi_j \geq \mu(x)\sum_{i=1}^n{\xi_i}^2$$
for all $n$-tuples of real numbers $(\xi_1 .. \xi_n)$. The operator
is said to be uniformly elliptic in a domain $\Omega$ if the
inequality holds for each point of $\Omega$ and if there is a
positive constant $\mu_0$ such that $\mu (x) \geq \mu_0$ for all $x$
in $\Omega$. Ellipticity of a more general second order operator is
defined via its second order term.

In the matrix language, the ellipticity condition asserts that the
symmetric matrix $[a^{ij}]$ is positive definite at each point $x$.

\beg{lem}[Hopf's strong maximum principle \cite{protter}]
\label{strongmax} Let $u$ satisfy the differential inequality
$$(L_c+h)u \geq 0 ~~with~~ h\leq 0$$ where $L_c$ is uniformly
elliptic in $\Omega$ and coefficients of $L_c$ and $h$ bounded. If
$u$ attains a nonnegative maximum at an interior point of $\Omega$,
then $u$ is constant.
\end{lem}

So if for example the maximum of $u$ is attained in the interior and
is $0$, then $u$ has to vanish. An application of this principle
provides us with a criterion of determining the sign of the Yamabe
Constant using Green's Functions:

 \begin{lem}[Green's Function Characterization for the Sign\cite{lom}]\label{greenchar}
 Let $(M,g)$ be a compact Riemannian $4$-manifold
with $Ker(\Delta + s/6)=0$, i.e. the Yamabe Laplacian has trivial
kernel, taking $\Delta=d^*d\cite{atgrn}$. Fix a point $y\in M$. Then
for the conformal class $[g]$ we have the following assertions : \\
1. It does not contain a metric of zero scalar curvature \\
2. It contains a metric of positive scalar curvature iff
$G_y(x)\neq0$ for all $x \in M-\{y\}$ \\
3. It contains a metric of negative scalar curvature iff $G_y(x)<0$
for some $x \in M-\{y\}$
\end{lem}
\begin{proof} Proceeding as in \cite{lom},
$[g]$ has three possibilities for its Yamabe Type, one of
$0$,$+$,$-$. Since the Yamabe Laplacian is conformally invariant as
acting on functions with conformal weight, we assume that either
$s=0$ or $s>0$ or else $s<0$ everywhere. \begin{itemize}   \item[]
\begin{itemize}
\item[$s=0$ :]Then $(\Delta+0/6)f=\Delta f=0$ is solved by any
nonzero constant function $f$. Therefore $Ker(\Delta+s/6)\neq 0$,
which is not our situation.
\item[$s>0$ :] For the smooth function $G_y : M-\{y\} \to
\mathbb{R}$~,~$G_y^{-1}((-\infty,a])$ is closed hence compact for
any $a \in \mathbb{R}$. Hence it has a minimum say at $m$ on
$M-\{y\}$. We also have $(\Delta+s/6)G_y=0$ on $M-\{y\}$. At the
minimum, choose normal coordinates so that $\Delta
G_y(m)=-\sum_{k=1}^4 \partial_k^2G_y(m)$. Second partial derivatives
are greater than or equal to zero, $\Delta G_y(m)\leq 0$ ~so~
$G_y(m)=-{6\over s}\Delta G_y(m)\geq 0$. We got nonnegativity, but
need positivity, so assume $G_y(m)=0$.\\ Then the maximum of $-G_y$
is attained and it is nonnegative with $(\Delta_c-s/6)(-G_y)=0\geq
0$. So the strong maximum principle(\ref{strongmax}) is applicable
and $-G_y\equiv 0$. This is impossible since $G_y(x)\to \infty$ as
$x\to y$, hence $m\neq 0$ and $G_y>0$. Note that the weak maximum
principle was not applicable since we had $G_y\geq 0$, implied
$\Delta_cG_y={s\over 6}G_y\geq 0$ though we got a minimum rather
than a maximum. Also note that $\nabla G_y(m)=0$ at a minimum though
this does not imply $div\nabla G_y(m)=0$.
\item[$s<0$ :]In this situation we have
$$\frac{1}{6}\int_M sG_ydV= \int_M(\Delta+s/6)G_ydV=\int_M \delta_ydV=1>0$$
implying $G_y<0$ at some point. Besides, at some other point it
should be zero since $G_y(x)\to +\infty$ as $x\to y$. \end{itemize}  \end{itemize}
\end{proof}

\section{Cohomological Characterization}\label{cohom} Now let $(M^4,g)$ be a
compact self-dual Riemannian manifold with the twistor space $Z$.
One of the basic facts of the twistor theory\cite{hitlin} is that
for any open set $U\subset M$ and the correponding inverse image
 $\wt U\subset Z$ in the twistor space, there is a natural isomorphism
$$pen : H^1 (\wt U, \mco(K^{1/2}))\stackrel{\sim}{\longrightarrow} \left\{ \beg{tabular}{c}smooth
complex-valued solutions\\ of~ $(\Delta  + s/6)u=0$ ~in~ $U$
\end{tabular} \right\}$$ which is called \emph{the Penrose
transform}\index{Penrose transform}\cite{bailsing,hitka,atgrn},
where $K=K_Z$. Since locally $\mco(K^{1/2}) \approx\mco(-2)$ e.g.
$Z=\mbb{CP}_3$, for a cohomology class $\psi \in H^1 (\wt U,
{\mathcal O}(K^{1/2}))$, the value of the corresponding function
$pen_\psi$ at $x\in U$ is obtained by restricting $\psi$ to the
twistor line $P_x\subset Z$ to obtain an element
$$pen_\psi(x)=\psi |_{P_x}\in H^1(P_x, {\mathcal O}(K^{1/2}))\approx  H^1(\bcp_1 , {\mathcal O}(-2) ) \approx \mbb C.$$
Note that  $pen_\psi$ is a section of a line bundle, but the choice
of a metric $g$ in the conformal class determines a canonical  trivialization
of  this line
bundle \cite{hitka}, and $pen_\psi$ then becomes an ordinary function.
Taking $U=M-\{y\}$ we have $(\Delta +s/6)G_y=0$ on $U$ in the uniquely presence of the conformal Green's functions(\ref{greatsmall}) and $G_y(x)$ is regarded as a function of $x$ corresponds to a canonical element $$pen^{-1}(G_y)\in H^1(Z-P_y, {\mathcal O}(K^{1/2}))$$
where $P_y$ is the twistor line over the point $y$.

What is this  interesting cohomology class?
The answer was discovered by Atiyah \cite{atgrn} involving
the \emph{Serre Class} of a complex submanifold. Which is a construction due to Serre \cite{serremod} and Horrocks \cite{horrocks}.
We now give the definition of the Serre class via the following lemma:


\begin{lem}[Serre-Horrocks Vector Bundle,Serre Class] \label{serhor}
Let $W$ be a (possibly non-compact) complex manifold, and let $V\subset W$ be a
closed complex submanifold of complex codimension $2$, and $N=N_{V/W}$ be the normal bundle of $V$.
For any holomorphic line bundle $L\to W$ satisfying
$$L|_V\approx \wedge^2 N ~~~\textnormal{and}~~~H^1(W, {\mathcal O}(L^*))=H^2(W, {\mathcal O}(L^*))=0$$
There is a rank-$2$ holomorphic vector bundle $E\to W$ called the
\emph{Serre-Horrocks bundle} of $(W,V,L)$ , together with a
holomorphic section $\zeta$ satisfying
$$\wedge^2 E \approx L~~~,~~~d\zeta|_V : N\stackrel{\sim}{\to} E~~~\textnormal{and}~~~\zeta =0 ~\textnormal{exactly on V}.$$
The pair $(E, \zeta)$ is unique up to isomorphism if we also impose
that the isomorphism $\det d\zeta: \wedge^2 N\to \wedge^2E|_V$
should agree with a given isomorphism $ \wedge^2 N\to L|_V$. They
also give rise to an extension
$$0\to {\mathcal O}(L^*)\to {\mathcal O}(E^*)
\stackrel{\cdot \zeta} {\to} {\mathscr I}_V\to 0,$$ the class of
which is defined to be the {\em Serre Class} $\lambda(V)\in
\tn{Ext}^{1}_W ({\mathscr I}_V, {\mathcal O}(L^*))$, where
${\mathscr I}_V$ is the ideal sheaf of $V$, and this extension
determines an element of $H^1(W-V, {\mathcal O}(L^*))$ by
restricting to $W-V$.
\end{lem}

\beg{proof} Consult \cite{lom} for a proof.
\end{proof}

For an alternative treatment of Serre's class via the Grothendieck class consult \cite{atgrn}.
We are now ready to state the answer of Atiyah :

\begin{thm}[Atiyah\cite{atgrn}] \label{atiyah}
Let $(M^4,g)$ be a compact self-dual Riemannian manifold
with twistor space $Z$, and assume that the conformally invariant
Laplace operator $\Delta=d^*d+s/6$ on $M$ has no global
nontrivial solution so that the Green's functions are well defined.
Let $y\in M$ be any point, and $P_y\subset Z$ be the corresponding twistor line.\\
 Then the image of the Serre class $\lambda (P_y)\in
 \tn{Ext}^{1}_Z ({\mathscr I}_{P_y}, {\mathcal O}(K^{1/2}))$
 in
 $H^1(Z-P_y, {\mathcal O}(K^{1/2}))$ is the Penrose transform of the Green's function $G_y$
 times a non-zero constant. More precisely
$$pen^{-1}(G_y)={1\over 4\pi^2}\lambda(P_y)$$
\end{thm}
Now thanks to this remarkable result of Atiyah, we can substitute
the Serre class for the Green's functions in our previous
characterization \ref{greenchar} and get rid of them to obtain a
better criterion for positivity as follows :
\begin{prop}[Cohomological Characterization , \cite{lom}]\label{nice}
Let $(M^4,g)$ be a compact self-dual Riemannian manifold with
twistor space $Z$. Let $P_y$ be a twistor line in $Z$. \\
Then the conformal class $[g]$ contains a metric of positive scalar curvature if and only if
  $H^1(Z, {\mathcal O}(K^{1/2}))=0$,  and the Serre-Horrocks vector bundle(\ref{serhor}) on $Z$ taking $L= K^{-1/2}$
associated to $P_y$ satisfies
$E|_{P_x}\approx {\mathcal O}(1)\oplus {\mathcal O}(1)$
for every twistor line $P_x$
\end{prop}
\beg{proof} $\Rightarrow$ : If a conformal class contains a metric
of positive scalar curvature $g$, then we can show that
$Ker(\Delta+{s\over 6})$ is trivial as follows: Let $(\Delta+{s\over
6})u=0$ for some smooth function $u:M\to\mbb R$ and $s>0$. Since $M$
is compact, $u$ has a minimum say at some point $m$. At the minimum
one has
$$\Delta u(m)=-\sum u_{kk}(m)\leq 0$$
because of the normal coordinates about $m$, modern Laplacian and
second derivative test. So that $$\Delta u=-{su\over 6}\leq 0
~~~\tn{implying}~~~ u\geq 0~~~ \tn{everywhere}.$$ If we integrate
over $M$ on gets $0$ for the Laplacian of a function so
$$0=\int_M\Delta u ~dV=\int_M -{su\over
6} dV$$ hence $$\int_M su~dV=0~~~\tn{implying}~~~u\equiv
0~~~\tn{since}~~~s>0$$ that is to say that the kernel is zero.

Remember the Penrose Transform map
$$pen : H^1(M, \mco(K^{1/2}))\stackrel{\sim}{\longrightarrow}
Ker(\Delta+{s\over 6})$$ implies that $H^1 (M, \mco(K^{1/2}))=0$,
also by Serre Duality $$H^2(M,K^{1/2})\approx
H^{0,2}_{\bar{\partial}}(M,K^{1/2})\stackrel{SD}{\approx}H^{3,1}_{\bar{\partial}}(M,K^{1/2*})^*\approx
H^1(M,K\otimes K^{-1/2})^*$$
$$=H^1(M,K^{1/2})^*=0$$ also
$$\wedge^2NP_y=\wedge^2\mco_{P_y}(1)\oplus\mco_{P_y}(1)=\bigoplus_{2=p+q}\wedge^p\mco(1)\otimes
\wedge^q\mco(1)=\wedge^1\mco(1)\otimes\wedge^1\mco(1)$$
$$=\mco_{\mbbp_1}(2)=K^{-1/2}|_{P_y}$$\\ since
$K^{-1/2}|_{P_y}=\mco_{\mbbp_3}(4)^{1/2}|_{P_y}=\mco_{P_y}(2).$ So
that the hypothesis for the Serre-Horrocks vector bundle
construction (\ref{serhor}) for $L=K^{-1/2}$ is satisfied. Then we
have the image of the Serre class
$$4\pi^2pen^{-1}(G_y)=\lambda(P_y)\in H^1(Z-P_y,K^{1/2})$$
So $$4\pi^2G_y(x)=pen_{\lambda(P_y)}(x)=\lambda(P_y)|_{P_x}\in
H^1(P_x,\mco(K^{1/2}))\approx \mbb C$$ where
$$H^1(P_x,\mco(K^{1/2}))\approx H^1(\mbb{CP}_1,\mco(-2))\approx H^0(\mbb{CP}_1,\Omega^1(\mco(-2)^*))
=H^0(\mbb{CP}_1,\mco)\approx \mbb C$$ By the Green's Function
Characterization (\ref{greenchar}) we know that $4\pi^2G_y(x)\neq
0$. So $\lambda(P_y)|_{P_x}\in H^1(P_x,\mco(K^{1/2}))$ is also
nonzero.

Since $\lambda(P_y)$ corresponds to the extension
$$0\to\mco(K^{1/2})\to\mco(E^*)\to \mc I_{P_y}\to 0$$
If we restrict to $Z-P_y$
$$0\to\mco(K^{1/2})\to\mco(E^*)\to \mco\to 0$$
dualizing we obtain
$$0\to\mco\to\mco(E)\to \mco(K^{-1/2})\to 0$$
now restricting this extension to $P_x$
$$0\to\mco_{\mbbp_1}\to\mco(E)|_{P_x}\to \mco(2)\to 0$$
So since $G_y(x)\neq 0$ , we expect that this extension is
nontrivial. Let's figure out the possibilities. First of all, by the
theorem of Grothendieck \cite{okonek}p22 every holomorphic vector
bundle over $\mbbp_1$ splits. In our case
$E|_{P_x}=\mco(k)\oplus\mco(l)$ for some $k,l\in\mbb Z$. Moreover if
we impose $k\geq l$, this splitting is uniquely
determined\cite{okonek}.

Secondly, any short exact sequence of vector bundles splits
topologically by \cite{okonek}p16. In our case, topologically we
have $E|_{P_x}\stackrel{t}{=}\mco\oplus\mco(2)$. So, setting the
chern classes to each other we have
$$c_1(E|_{P_x})[P_x]=c_1(\mco(k)\oplus\mco(l))[\mbbp_1]=c_1\mco(k)+c_1\mco(l)[\mbbp_1]=k+l$$
equal to
$$c_1(E|_{P_x})[P_x]=c_1(\mco\oplus\mco(2))[\mbbp_1]=c_1\mco+c_1\mco(2)[\mbbp_1]=0+2=2.$$
Hence $l=2-k$. We now have $E|_{P_x}=\mco(k)\oplus\mco(2-k)$. Our
extension becomes $$0\to\mco_{\mbbp_1}\to\mco(k)\oplus\mco(2-k)\to
\mco(2)\to 0$$ The inclusion
$\mco\hookrightarrow\mco(k)\oplus\mco(2-k)$ gives a trivial
holomorphic subbundle. It has one complex dimensional space of
sections. So these sections are automatically sections of
$\mco(k)\oplus\mco(2-k)$, too. This implies
$$0\neq H^0(\mco(k)\oplus\mco(2-k))=H^0(\mco(k))\oplus H^0(\mco(2-k))$$
Imposing $k,2-k\geq 0$ by the Kodaira Vanishing Theorem\cite{gh}
since the direct sum elements $\mco(k)$ and $\mco(2-k)$ should
possess sections. Also, from uniqueness $k\geq l=2-k$. Altogether we
have
 $2\geq k\geq 1$. From the two choices, $k=2$ gives the trivial
 extension $\mco(2)\oplus\mco$, $k=1$ gives the nontrivial extension
 $E|_{P_x}=\mco(1)\oplus\mco(1)$ as we expected. See the following remark for existence.

$\Leftarrow$ : For the converse, if $E|_{P_x}=\mco(1)\oplus\mco(1)$
then we already showed that this is the nontrivial extension hence
$G_y(x)\neq 0$, so that the scalar curvature is positive by the
Green's Function Characterization (\ref{greenchar}) \end{proof}

\beg{rmk}The nontrivial extension of $\mco$ by $\mco(2)$ exists
 by the Euler exact sequence\index{Euler sequence}
$$0\to \mco\to\mco(1)^{\oplus n+1}\stackrel{\mc E}{\to} T'\mbbp^n\to 0$$
\cite{gh}p409 for $n=1$. Alternatively, the maps $i:\rho\mapsto
(\rho Z_0,\rho Z_1)$ and $j:(u,v)\mapsto uZ_1-vZ_0$ for coordinates
$[Z_0:Z_1]$ on $\mbbp_1$
 yields the exact sheaf sequence
$$0\to\mco(-1)\stackrel{i}{\to}\mco\oplus\mco\stackrel{j}{\to}\mco(1)\to
0$$ tensoring with $\mco(1)$ produces the nontrivial
$\mco(1)\oplus\mco(1)$ extension. Since we have a unique nontrivial
extension, this shows
$$\tn{Ext}^1(\mco(2),\mco)=\mbb C$$
used in $\cite{atgrn}$ to classify the extensions. On the other hand
$$H^1(Hom(\mco(2),\mco))=H^1(\mco(2)^*\otimes\mco)=H^1(\mco_{\mbbp_1}(-2))$$ $$=H^0(\mco(2)\otimes\mco(-2)^*)=
H^0(\mbbp_1,\mco)=\mbb C$$ used in $\cite{DF}$ to classify the
extensions. So, our computation verifies the isomorphism
$$\tn{Ext}^q(M,\mc F,\mc G)\approx H^q(M,\mc F^*\otimes_\mco \mc
G)$$ for locally free sheaves or vector
bundles 
for $q=1$. See \cite{gh}$p706$.\\
Here, \tn{Ext} stands for what is called the {\em global
Ext}\index{Ext, global} group usually defined to be the
hypercohomology of the complex of sheaves associated to a global
syzygy for $\mc F$. Though practically usually computed via the
spectral sequence to be
$$\tn{Ext}^k(\mc F,\mc G)=H^0(\underline{Ext}^k_\mco(\mc F,\mc G))$$
under some vanishing conditions\cite{gh}. \hfill $\square$
\end{rmk}

\section{The Sign of the Scalar Curvature} \label{scalar} We are now
ready to approach  the problem of determining  the sign of the
Yamabe constant for the self-dual conformal classes constructed in
Theorem (\ref{thmb}). The techniques used here are analogous to the
ones used by LeBrun in \cite{lom}.

\begin{thm}\label{pos}
Let $(M_1,g_1)$ and $(M_2,g_2)$be compact self-dual Riemannian
 $4$-manifolds with $H^2(Z_i, {\mathcal O}(TZ_i))=0$ for their twistor spaces.
Moreover suppose that they have positive scalar curvature.

Then,  for all sufficiently small ${\mathfrak t}>0$, the self-dual
conformal class $[g_{\mathfrak t}]$ obtained on $M_1 \# M_2$ by the
Donaldson-Friedman Theorem (\ref{thmb}) contains a metric of
positive scalar curvature.
\end{thm}

\beg{proof} Pick a point $y \in (M_1 \# M_2)\backslash M_1$.
Consider the real twistor line $P_y\subset \wtz_2$, and extend this
as a $1$-parameter family of twistor lines  in $P_{y_{\mathfrak
t}}\subset Z_{\mathfrak t}$ for ${\mathfrak t}$ near $0\in \mathbb
C$ and such that
 $P_{y_{\mathfrak t}}$ is a real twistor line
for ${\mathfrak t}$ real. By shrinking ${\mathcal U}$ if needed, we
may arrange that ${\mathcal P}= \cup_{\mathfrak t} P_{y_{\mathfrak
t}}$ is a closed codimension-$2$ submanifold of ${\mathcal Z}$ and
$H^1 ({\mathcal Z} , {\mathcal O}(L^*))=H^2 ({\mathcal Z} ,
{\mathcal O}(L^*))=0$ by the Vanishing Theorem
 (\ref{vanishing}). Next we check that $L|_\mathcal
 P\approx \wedge^2N_{\mathcal
 P}$. Over a twistor line $P_{y_\mathfrak t}$ we have $$\wedge^2N_\mathcal P|_{P_{y_\mathfrak t}}
 =\wedge^2 (\mco(1)\oplus\mco(1))=\mco_{P_{y_\mathfrak t}}(2)$$ by considering the first
 Chern classes. On the other hand, notice that the restriction of $L^*$ to any smooth fiber
$Z_{\mathfrak t}$, ${\mathfrak t}\neq 0$ is simply $K^{1/2}$ :
$$L^*|_{Z_\mathfrak t}=({1\over 2}K_\mcz-\wtz_1)|_{Z_\mathfrak t}={1\over 2}K_\mcz|_{Z_\mathfrak t}=
{1\over 2}(K_{Z_\mathfrak t}-Z_\mathfrak t)|_{Z_\mathfrak t}=
{1\over 2}K_{Z_\mathfrak t}|_{Z_\mathfrak t}.$$ Here,
$\wtz_1|_{Z_\mathfrak t}=0$ because of the fact that $\wtz_1$ and
${Z_\mathfrak t}$ does not intersect for ${\mathfrak t}\neq 0$. The
normal bundle of $Z_\mathfrak t$ is trivial, because of the fact
that we have a standard deformation. Then
$$L|_{P_{y_\mathfrak t}}=
K^{-{1/2}}_{Z_\mathfrak t}|_{P_{y_\mathfrak t}}=TF|_{P_{y_\mathfrak
t}}=\mco_{P_{y_\mathfrak t}}(2) ~~~\textnormal{for}~~~ \mathfrak
t\neq 0$$ since ~$TF$ of Sec (\ref{natural}) is the square-root of
the anti-canonical bundle. For the case $\mathfrak t=0$, we need the
fact that ~$L^*|_{\wt{Z}_2}=\pi^*K^{1/2}_{Z_2}$~ which we have
computed in the step \ref{hitvan} of the proof of the vanishing
theorem (\ref{vanishing}). This yields
$$L|_{P_{y_0}}=\pi^*K^{-1/2}_{Z_2}|_{\wtz_2}|_{P_{y_0}}=\mco_{P_{y_0}}(2).$$

Then the Serre-Horrocks construction (\ref{serhor}) is available to
obtain the holomorphic vector bundle $E\to {\mathcal Z}$  and a
holomorphic section $\zeta$ vanishing exactly along $\mathcal P$,
also, the corresponding extension $$0\to {\mathcal O}(L^*)\to
{\mathcal O}(E^*)\stackrel{\cdot\zeta}{\to} {\mathscr I}_{\mathcal
P}\to 0$$ gives us the Serre class $\lambda({\mathcal P})\in
H^1({\mathcal Z}-{\mathcal P}, {\mathcal O}(L^*))$.

Since $L^*|_{Z_\mathfrak t}=K^{1/2}_{Z_\mathfrak t}$ for $\mathfrak
t\neq 0$ by the above computation, Proposition (\ref{atiyah}) of
Atiyah tells us that the restriction of $\lambda ({\mathcal P})$ to
$Z_{\mathfrak t}$, ${\mathfrak t} > 0$, has Penrose transform equal
to a positive constant times the conformal Green's function of $(
M_1\# M_2, g_{\mathfrak t}, y_{\mathfrak t})$ for any ${\mathfrak t}
> 0$.

Now, we will restrict $(E,\zeta )$ to the two components of the
divisor $Z_0$. We begin by restricting to $\wt{Z}_2$. We have
$L|_{P_{y_0}}=\mco_{P_0}(2)=\wedge^2N{P_{y_0}}$ and
$$H^k(\wtz_2,L^*)=H^k(\wtz_2,\pi^*K^{1/2}_{Z_2})=H^k(Z_2,\pi_*\pi^*K^{1/2}_{Z_2})
=H^k(Z_2,K^{1/2}_{Z_2})=0$$ for $k=1,2$ because of the projection
lemma, Leray spectral sequence and the Hitchin's Vanishing theorem
for positive scalar curvature on $M_2$. So that we have the
Serre-Horrocks bundle for the triple $(\wtz_2 , P_{y_0} ,
L|_{\wtz_2}=\pi^*K^{-1/2}_{Z_2})$. On the other hand it is possible
to construct the Serre-Horrocks bundle $E_2$ for the triple
$(Z_2,P_{y_0},K_{Z_2}^{-1/2})$ for which all conditions are already
checked to be satisfied. In the construction of these Serre-Horrocks
bundles, if we stick to a chosen isomorphism $\wedge^2N \to
L|_{P_{y_0}}$, these bundles are going to be isomorphic by
(\ref{serhor}). 
The splitting type of $E$ on the twistor lines corresponding to the
points in $M_2-\{y_0,p_2\}$ supposed to be the same as the splitting
type of $E_2$, which is $\mco(1)\oplus\mco(1)$ since $Z_2$ already
admits a self-dual metric of positive scalar curvature.

Secondly, we restrict $(E,\zeta)$ to $\wtz_1$. Alternatively we
restrict the Serre class $\lambda(\mc P)$ to $H^1(\wtz_1,\mco(L^*))$
where \bct $L^*|_{\wtz_1}=\frac{1}{2}K_\mcz-\wtz_1|_{\wtz_1}
=\frac{1}{2}K_\mcz+Q|_{\wtz_1}\stackrel{adj}=\frac{1}{2}(K_{\wtz_1}-\wtz_1)+Q|_{\wtz_1}
=\frac{1}{2}(K_{\wtz_1}+Q)+Q|_{\wtz_1}=\frac{1}{2}(\pi^*K_{Z_1}+2Q)+Q|_{\wtz_1}=\pi^*{1\over
2}K_{Z_1}+2Q|_{\wtz_1},$ \ect and show that it is non-zero on every
real twistor line away from $Q$ here. Remember that we have the the
restriction isomorphism obtained in the step \ref{restrictiontoq} of
the proof of the vanishing theorem (\ref{vanishing})
$$H^1(\mco_{\wtz_1}(L^*))\stackrel{\sim}{\longrightarrow}H^1(\mco_Q(L^*))\approx \mbb C$$
as a consequence of Hitchin's Vanishing theorems for positive scalar
curvature on $M_1$, as mentioned in the step \ref{tech}, and
$H^1(\mco_Q(L^*))=H^1({\mbb P}_1 \times {\mbb P}_1,\mco(-2,0))=\mbb
C$, as computed in the step \ref{restrto quadric cohomology}. This
shows that if there is a rational curve of $Q$ on which the Serre
class is non-zero, then this class is non-zero and a generator of
$H^1(\mco_{\wtz_1}(L^*))$. The Serre-Horrocks bundle construction on
$Z_2$ shows us that $E|_{C_2}= \mco(1)\oplus\mco(1)$ where $C_2$ is
the twistor line on which the blow up is done. We know that
$Q=\mbbp_1\times\mbbp_1\approx \mbb P(NC_2)$. So that the
exceptional divisor has one set of rational curves which are the
fibers, and another set of rational curves, coming from the sections
of the projective bundle $\mbbp (NC_2)$. Take the zero section of
$\mbbp(NC_2)$, on which $E$ has a splitting type
$\mco(1)\oplus\mco(1)$. So over the zero section in $Q$, $E$ is
going to be the same, hence non-trivial splitting type. This shows
that over this rational curve on $Q$, the Serre-class is nonzero.
Hence by the isomorphism above, the Serre-class is the (up to
constant) nontrivial class in $H^1(\wtz_1,\mco(L^*))\approx\mbb C$.

Next we have to show that this non-trivial class is non-zero on
every real twistor line in $\wtz_1-Q$ or $Z_1-C_1$\footnote{Thanks
to C.LeBrun for this idea.}. For this purpose consider the
Serre-Horrocks vector bundle $E_1$ and its section $\zeta_1$ for the
triple $(Z_1,C_1,K_{Z_1}^{-1/2})$, so that $\pi^*\zeta_1$ is a
section of $\pi^*E_1$ vanishing exactly along $Q$. Remember the
construction of the line bundle associated to the divisor $Q$ in
$\wtz_1$ \cite{gh}. Consider the local defining functions $s_\alpha
\in \mf M^*(U_\alpha)$\footnote{Here, $\mf M^*$ stands for the
multiplicative sheaf of meromorphic functions which are not
identically zero, in the convention of \cite{gh}. Actually the local
defining functions here are holomorphic because Q is effective.} of
$Q$ over some open cover $\{U_\alpha\}$ of $\wtz_1$. These functions
are holomorphic and vanish to first order along $Q$. Then the
corresponding line bundle is constructed via the transition
functions ~$g_{\alpha\beta}=s_\alpha\ / s_\beta$. Since $s_\alpha$'s
transform according to the transition functions, they constitute a
holomorphic section $s$ of this line bundle $[Q]$, which vanish up
to first order along $Q$. Local holomorphic sections of this bundle
is denoted by $\mco([Q])$ and they are local functions with simple
poles along $Q$. If we multiply $\pi^*\zeta_1$ with these functions,
we will get a holomorphic section of $\pi^*E_1$ on the corresponding
local open set, since $\zeta_1$ has a non-degenerate zero on $Q$, so
that it vanishes up to degree $1$, there. This guarantees that the
map is one to one, and the multiplication embeds $\mco([Q])$ into
$\pi^*E_1$. The quotient has rank $1$, and the transition functions
of $\pi^*E_1$ relative to a suitable trivialization will then look
like
$$\left(\begin{array}{cc}                       g_{\alpha\beta} & k_{\alpha\beta} \\
                                                0 &
                                                d_{\alpha\beta}\cdot
                                                g_{\alpha\beta}^{-1}\end{array}\right)$$
where $d_{\alpha\beta}$ stands for the determinant of the transition
matrix of the bundle $\pi^*E_1$ in this coordinate chart. Since the
bundle ~$\det\pi^*E_1\otimes[Q]^{-1}$ has the right transition
functions, it is isomorphic to the quotient bundle, hence we have
the following exact sequence $$0\to [Q] \to \pi^* E_1 \to
\pi^*K^{-1/2}  \otimes [Q]^{-1} \to 0$$ since ~$\det
E_1=K^{-1/2}_{Z_1}$ as an essential feature of the Serre-Horrocks
construction. This extension of line bundles is classified by an
element in
$$\mrm{Ext}^1_{\wtz_1}(\pi^*K^{-1/2}  \otimes [Q]^{-1},[Q])\approx
H^1(\wtz_1,\pi^*K^{1/2}  \otimes [Q]^{2})$$ by \cite{atgrn}. If we
restrict our exact sequence to ~$\wtz_1-Q=Z_1-C_1$, since the bundle
$[Q]$ is trivial on the complement of $Q$, this extension class will
be the Serre class of the triple $(Z_1,C_1,K_{Z_1}^{-1/2})$.
Finally, since $M_1$ has positive scalar curvature, this class is
nonzero on every real twistor line in $Z_1-C_1$. So that
non-triviality of the class forced non-triviality over the real
twistor lines. In other words $E$ has a non-trivial splitting type
over the real twistor lines of $\wtz_1$.





So we showed that the Serre-Horrocks vector bundle $E$ determined by
$\lambda ({\mathcal P})$ splits as ${\mathcal O}(1)\oplus {\mathcal
O}(1)$ on all  the $\sigma_0$-invariant rational curves in $Z_0$
which are limits of real twistor lines in ${\mathcal Z}_{\mathfrak
t}$ as ${\mathfrak t}\to 0$. It therefore has the same splitting
type on all the real twistor lines of ${\mathcal Z}_{\mathfrak t}$
for ${\mathfrak t}$ small. Besides, $$h^j(Z_t,\mco(L^*)) \leq
h^j(Z_0,\mco(L^*))=0 ~~~\textnormal{for} ~~~j=1,2$$ by the
semi-continuity principle and the proof of the vanishing theorem
(\ref{vanishing}). So that via~ $L^*|_{Z_\mathfrak t}\approx
K^{1/2}$, $$H^1(Z_{\mathfrak t},\mco(K^{1/2}))\approx
Ker(\Delta+{s\over 6})=0.$$

Since the two conditions are satisfied, Cohomological
characterization (\ref{nice})
guarantees the positivity of the conformal class.
\end{proof}

\bigskip

\vspace{5mm}
\small \beg{flushleft} \textsc{Department of Mathematics, State University of New York, Stony Brook}\\
\textit{E-mail address :} \texttt{\textbf{kalafat@math.sunysb.edu}}
\end{flushleft}


\end{document}